\documentclass[reqno,9pt]{amsart}
\usepackage[applemac]{inputenc}

\usepackage{amsmath}
\usepackage{amssymb}
\usepackage{amsfonts}
\usepackage{amsthm}
\usepackage{multirow}
\usepackage{enumerate}
\usepackage{color}
\usepackage{dsfont}

\usepackage{setspace}
\newcommand{\R}{\mathbb{R}}
\newcommand{\C}{\mathbb{C}}
\newcommand{\f}{\rightarrow}
\newcommand{\deb}{\bar\partial}
\newcommand{\de}{\partial}

\newcommand{\lmb}{\lambda}
\newcommand{\ov}[1]{\overline{#1}}
\newcommand{\w}[1]{\widetilde{#1}}

\newcommand{\M}{\mathcal M}

\newcommand{\hyp}{\operatorname{hyp}}

\newcommand{\tr}{\operatorname{tr}}

\newcommand{\arctanh}{\operatorname{arctanh}}

\newcommand{\Vol}{\operatorname{Vol}}

\newcommand{\rk}{\operatorname{rk}}

\newcommand{\Entv}{\operatorname{Ent}}

\newcommand{\ep}{{\varepsilon}}

\newcommand{\FF}{{\mathcal F}}

\newcommand{\reg}{{\operatorname{reg}}}
\newcommand{\di}{{\operatorname{d}}}
\newcommand{\om}{{\mathrm{hyp}}}

\newtheorem{prop}{Proposition}
\newtheorem{thm}{Theorem}

\newtheorem{ex}{Example}
\newtheorem{rmk}{Remark}

\begin{document}

\author[R. Mossa]{Roberto Mossa}

\address{Laboratoire de MathŽmatiques Jean Leray (UMR 6629) CNRS, 2 rue de la Houssinire -- B.P. 92208 -- F-44322 Nantes Cedex 3, France}

\email{roberto.mossa@gmail.com}

\title[The volume entropy of local HSSNCT]
{The volume entropy of local Hermitian symmetric space of noncompact type}

\begin{abstract}
We calculate the volume entropy of local Hermitian symmetric spaces of noncompact type in terms of its invariant $r$, $a$, $b$.
\end{abstract}

\maketitle

\section{introduction and statement of the main result}

Let $( M,g)$ be a compact Riemennian manifold with Riemannian universal covering $(\w M, \w g)$. 
The \emph{volume entropy} of $M$ also called the \emph{exponential rate of the volume growth}, is the asymptotic invariant defined by 
\begin{equation}\label{def vol ent}
\begin{split}
\Entv( M, g)=\lim_{t \f \infty} \frac{1}{t}\log\Vol\left( B_p\left(t\right)\right),
\end{split}
\end{equation}
here $\Vol\left( B_p\left(t\right)\right)$ is the volume of the geodesic ball $B_p(t)\subset \w M$ of center $p$ and radius $t$. Due to the compactness this limit exists, moreover it does not depend on the point $p\in \w M$ (see \cite{manning}).
Integrating by parts \eqref{def vol ent}, the entropy can be equivalently defined  as the infimum of the positive constants $c$ such that the integral 
\begin{equation*}\label{def vol ent int00}
\int_{ \w M} e^{-c\, \di(p,\cdot)}\, d\tilde v
\end{equation*}
converges, where $\di(p,\cdot)$ is the geodesic distance from $p$ and  $d\tilde v$ is the volume form of $\w M$ induced by $ \w g$,
namely 
\begin{equation}\label{def vol ent int}
\Entv( M, g)=\left\{\inf_{c\in \R}\mid \int_{ \w M} e^{-c\, \di(p,\cdot)}\, d\tilde v<\infty\right\}.
\end{equation}
 Notice that the volume  entropy $\Entv(M, g)$ only depends on the Riemannian universal covering ,
  therefore it make sense to define the volume entropy of $(\w M, \w g)$ as that of $(M, g)$, i.e.
\begin{equation*}\label{def vol ent covering}
\begin{split}
\Entv( \w M, \w g):=\Entv( M,g).
\end{split}
\end{equation*}

The study of the volume entropy starts with A. Manning in \cite{manning} who  shows that it is always bounded above by the topological entropy of the geodesic flow on $(M,g)$. Moreover, if $g$ has nonpositive sectional curvature then the volume entropy and the topological entropy coincide. The  volume entropy $\Entv(M,g)$ carries many geometric  informations about $(M, g)$.  For instance G. Besson, G. Courtois and S. Gallot in their celebrated paper \cite{bcg1} show that compact locally symmetric spaces of negative curvature are characterized by the minimality of the volume entropy among all homotopy equivalent Riemannian manifolds with the same volume. Further,  they provide another proof of Mostow's rigidity theorem in the case of rank 1 symmetric spaces (the reader is referred  to \cite{bcg1}, \cite{bcg2} and \cite{bcg3} for details).

The aim of this paper is to compute the entropy of compact quotients of irreducible Hermitian symmetric spaces of noncompact type (from now on HSSNCT) in terms of their invariants  $r$, $a$ and $b$ (see next section).
\begin{thm}\label{alternative}
Let $(\Omega,g_{\om})$ be an irreducible HSSNCT endowed with the hyperbolic metric $g_{\om}$ of holomorphic sectional curvature between $0$ and $-4$. Then its volume entropy is given by 
\begin{equation}\label{entv hssnct}
\operatorname{Ent}\left(\Omega, g_{\om}\right)=2\,\sqrt{\sum _{j=1}^{r}
 \left( b+1+a \left( r-j \right)  \right) ^{2}},
\end{equation}
where $r,a$ and $b$ are the invariant associated to $\Omega$.
\end{thm}
\noindent
As the entropy $\operatorname{Ent}\left(\Omega, g_{\om}\right)$ of an irreducible HSSNCT is determined by $r$, $a$ and $b$, it is natural to ask if it determines the domain $\Omega$, in general that is not true, indeed  the domains $\Omega_I[12,2]$ and $\Omega_{IV}[18]$ (whose corresponding constant are respectively $(r=2,\,a=2,\,b=10)$ and $(r=2,\,a=16,\,b=0)$) have the same volume entropy. Nevertheless we believe that once fixed the dimension $d$ the volume entropy $\operatorname{Ent}\left(\Omega, g_{\om}\right)$ determines the irreducible HSSNCT $\Omega$.

 Let now $\Omega=\Omega_1 \times \dots \times \Omega_\ell=\prod_{k=1}^\ell \Omega_k$ be a product of irreducible bounded symmetric domains. The above theorem extends to the following
\begin{thm}\label{thm prod}
The volume entropy of $\Omega=\prod_{k=1}^\ell \Omega_k$ is given by
\[
\Entv(\Omega,g)=\sqrt{\sum_{k=1}^\ell\Entv^2(\Omega_k,g_k)}
\]
where $g_k$ is the hyperbolic metric associated to $\Omega_k$.
\end{thm}

The proofs of Theorem \ref{alternative} and \ref{thm prod} use the  properties of the \emph{symplectic duality map} introduced by A. Loi and A. Di Scala in \cite{Loi} (see also \cite{exponential} and \cite{unicdual}) and use the theory of Jordan triple system associated to an HSSNCT (see section below).
We finally point out that in the Appendix of \cite{kernel}  one can find the computation of the first eigenvalue  (and hence the value of the volume entropy) of HSSNCT in terms of its root systems. In principle this could provide an alternative proof of Theorem  \ref{thm prod}  once one is able to pass from the root systems description to the one in terms of the invariants  $r, a, b$;  this seems not to be an easy task.
\\

The paper consists of three other sections. In the next one we recall the basic material about Hermitian positive Jordan triple systems, Hermitian symmetric spaces of noncompact type and the description of the geodesic distance in polar coordinates of these spaces. Section \ref{sect 3} and \ref{sec 4} are dedicated to the proof of Theorem \ref{alternative} and \ref{thm prod} respectively.\\

\noindent {\bf Acknowledgments}. 
The author would like to thank Professor Andrea Loi for various stimulating discussions, Professor Gilles Carron for his interest in my work and his comments, Professor Sylvestre Gallot for his help and continuous encouragement  and Professor Guy Roos for his interest in my work and for pointing me out the volume form of an Hermitian symmetric space of noncompact type in terms of Jordan triple systems.

%%%%%%%%%%%%%%%%%%%%%%%%%%%%%%%%%%%%%%%%%%%%%%%

\section{Hermitian symmetric spaces of noncompact type and Hermitian positive Jordan triple systems}
We refer the reader  to \cite{roos} (see also \cite{loos}) for more details on Hermitian symmetric spaces of noncompact type and Hermitian positive Jordan triple systems (from now on HPJTS).

\subsection{Definitions and notations}
An Hermitian Jordan triple system is a  pair $\left({\mathcal M},
\{ ,  ,\}\right)$, where ${\mathcal M}$ is a complex vector space and $\{
,  ,\}$ is a map
\[
\{ ,  ,\}:{\mathcal M}\times {\mathcal M}\times {\mathcal M} \rightarrow {\mathcal M}
\]
\[
\left(u, v, w\right)\mapsto \{u, v, w\}
\]
which is ${\C}$-bilinear and symmetric in $u$ and $w$, ${\C}$-antilinear in $v$ and such that the following \emph{ Jordan identity} holds:
\[
\{x, y, \{u, v, w\}\}-\{u, v, \{x, y, w\}\}= \{\{x, y, u\}, v,
w\}-\{u, \{v, x, y\}, w\}.
\]
For $x,y,z \in \M$ consider  the following operators
\[
T\left(x,y\right)z =\left\{  x,y,z\right\} 
\]
\[
Q\left(x,z\right)y =\left\{  x,y,z\right\}  
\]
\[
Q\left(x,x\right) =2\,Q\left(x\right)\label{D3}\\
\]
\[
B\left(x,y\right) =\operatorname{id}_{\mathcal M}-T\left(x,y\right)+Q\left(x\right)Q\left(y\right). \label{D4}
\]
The operators $B\left(x,y\right)$ and $T\left(x,y\right)$ are $\mathbb{C}$-linear, the operator
$Q\left(x\right)$ is $\mathbb{C}$-antilinear. $B\left(x,y\right)$ is called the \emph {Bergman operator}.
For $z\in \M$, the \emph{odd powers} $z^{\left(2p+1\right)}$ of $z$ in the Jordan triple
system $\M$ are defined by
\[
z^{\left(1\right)}=z \qquad z^{\left(2p+1\right)}=Q\left(z\right)z^{\left(2p-1\right)}. \label{D6}
\]
An Hermitian Jordan triple system is called  \emph {positive} if the Hermitian form
\[
\left(  u\mid v\right)  =\tr T\left(u,v\right) \label{D5}
\]
is positive definite. An element $c \in \M$ is called \emph {tripotent} if
$\{c,c,c\}=2c$. Two tripotents $c_1$ and $c_2$ are called \emph {(strongly)
orthogonal} if $T\left(c_1, c_2\right)=0$.

\subsection{HSSNT associated to HPJTS}
An Hermitian symmetric space of noncompact type $\Omega$ is uniquely determined by a triple of integers $(r,a,b)$, where $r$ represents the rank of $\Omega$ and $a$ and $b$ are positive integers. The dimension $d$ of $\Omega$ satisfies $2\,d=r\left( 2\, b + 2 + a\left(r-1\right) \right)$ and the genus $\gamma$ of $\Omega$ is given by $\gamma=\left(r-1\right)a+b+2$. Observe that $(\Omega,g_{\om})=\C H^n$ if and only if its rank is equal to $1$.
The table below summarizes the numerical invariants and the dimension of irreducible HSSNCT according to its type (for a more detailed description of this invariants, which is not necessary in our approach, see e.g. \cite{arazy}, \cite{zhang}).

\begin{tabular}{|p{0.1\textwidth}|p{0.11\textwidth}|p{0.11\textwidth}|p{0.1\textwidth}|p{0.1\textwidth}|p{0.1\textwidth}|p{0.1\textwidth}|}
\hline
\small Type         & $\Omega_{I}[{n,m}]$  &  $\Omega_{II}[{n}]$ & $\Omega_{III}[{n}]$ & $\Omega_{IV}[{n}]$ & $\Omega_{V}$ & $\Omega_{VI}$\\

                 &\small  $\{n\leq m\}$              &\small  $\{5 \leq n\}$           &\small $\{2 \leq n\}$ &\small $\{5 \leq n\}$                 &\small&\small\\ \hline
$d$ & $nm$ & $\frac{(n-1)n}{2}$               & $\frac{(n+1)n}{2}$ & $n$ & $16$ & $27$ \\ \hline
$r$ & $n$     & $\left[ \frac{n}{2} \right]$     & $n$                        & $2$ & $2$   & $3$   \\ \hline
$a$ & \small{$2$, if $2\leq n$ $0$, if $1=n$} & $4$     & $1$              & $n-2$ & $6$   & $8$   \\ \hline
$b$ & $m-n$  & \small{$0$, if $n$ even $2$, if $n$ odd} & $0$  & $0$ & $4$   & $0$   \\ \hline
$\gamma$ & $m+n$  & $2n-2$ & $n+1$  & $n$ & $12$   & $18$   \\ \hline
\end{tabular}\\

%\newpage
M. Koecher (\cite{Koecher1}, \cite{Koecher2}) discovered that to every HPJTS
$\left(\M, \{ ,  ,\}\right)$ one can associate an Hermitian symmetric
space of noncompact type, i.e. a bounded symmetric domain $\Omega$
centered at the origin $0\in \M$. The domain $\Omega$ is defined as the connected component containing the origin of   the set of all $u\in {\M}$ such that $B\left(u, u\right)$ is positive definite with respect to the Hermitian form $\left(  u\mid v\right)  =\tr T\left(u,v\right) \label{D5}$. \emph{We will always consider such a domain in its (unique up to linear isomorphism) circled realization.} Suppose that $\M$ is \emph{simple} (i.e. $\Omega$ is irreducible).
 The \emph{flat} form $\omega_0$ is defined by
\[
\omega_0= -\frac{i}{2\,\gamma} \partial \bar \partial \log\tr T(x,\,x).
\]
If $\left(  z_{1},\ldots,z_{d}\right)  $ ($d=\dim( \M)$) are orthonormal
coordinates for the Hermitian product $\left(  u\mid v\right)  $, then%
\[
\omega_{0}=\frac{i}{2}\sum_{m=1}^{d}{d}z_{m}\wedge
{d}\overline{z}_{m}.
\]
The reproducing kernel ${K_\Omega}$ of $\Omega$, with respect $\omega_0$ is given by
\begin{equation*}\label{KOB}
\left(K_\Omega\left(z,\, \bar z\right)\right)^{-1}=C \det B\left(z,\,z\right),
\end{equation*}
where $C=\int_\Omega \frac{\omega_0^n}{n!}$.
When $\Omega$ is irreducible
\[
\omega_{\om}= -\frac{i}{2\,\gamma} \partial \bar \partial \log\det B= \frac{\omega_B}{\gamma},
\]
is the \emph{hyperbolic} form on $\Omega$, with associated hyperbolic metric $g_{\om}$ of holomorphic sectional curvature between $0$ and $-4$ and $\omega_B=-\frac{i}{2} \partial \bar \partial \log\det B$ is the Bergman form associated to the Bergman metric on $\Omega$. The HSSNCT associated to $\M$ is $(\Omega, g_\om)$.

\begin{ex}The hyperbolic space. \rm
Suppose $\M=\C$ and $\{x,\,y,\,z\}=x\,\ov y\, z$. The associated operators are 
\[T(x,\,y)=x\,\ov v\]
\[Q(x)\,y=x^2\,\ov y\]
\[B(x,\,y)=(1-x\,\ov y)^2.\]
Therefore the domain \[\Omega=\{z\in \C \,:\, B(z,\,z)>0\}\] is the unit disc of $\C$ and the hyperbolic form is 
\begin{equation}\label{hypform}
\omega_{\om}=-\frac{i}{2}\de \deb \log (1-|z|^2).
\end{equation} 
Namely $(\Omega, g_{\om})$ is the hyperbolic space $\C H^1$ with the hyperbolic metric of curvature $-4$.
\end{ex}

The HPJTS $\left({\M}, \{ ,  ,\}\right)$ can be recovered by its
associated HSSNT $\Omega$ by defining ${\M}=T_0 \Omega$ (the tangent space to the origin of $\Omega$) and
\begin{equation*}\label{trcurv}
\{u, v, w\}=-\frac{1}{2}\left(R_0\left(u, v\right)w+J_0\,R_0\left(u, J_0\,v\right)w\right),
\end{equation*}
where $R_0$ (resp. $J_0$) is the curvature tensor of the Bergman metric (resp. the complex structure) of $\Omega$ evaluated at the origin.
The reader is referred   to  Proposition III.2.7 in  \cite{Bertram} for the proof of (\ref{trcurv}). For more informations on the correspondence between
HPJTS and HSSNT we  refer also  to p. 85 in Satake's book \cite{satake}.

\subsection{Totally geodesic submanifolds of HSSNT}
We have the following result:
\begin{prop} \label{sub}
Let ${\Omega}$ be a HSSNT and let ${\M}$ be its associated HPJTS. Then there exists a one to one
correspondence between (complete) complex totally geodesic
submanifolds through the origin and sub-HPJTS of ${\M}$. This correspondence
sends $T \subset {\Omega}$ to  ${\mathcal T} \subset {\M}$, where
${\mathcal T}$ denotes the HPJTS associated to $T$.
\end{prop}

\subsection{Spectral decomposition and Functional calculus}
Let $\M$ be a HPJTS. Each element $z\in \M$ has a unique \emph{spectral decomposition}
\[
z=\lambda_{1}\,c_{1}+\cdots+\lambda_{s}\,c_{s}\qquad\left(0<\lambda_{1}<\cdots
<\lambda_{s}\right), \label{D7}
\]
where $\left(c_{1},\ldots,\,c_{s}\right)$ is a sequence of pairwise
orthogonal tripotents and the $\lambda_j$ are real number called eigenvalues of $z$. The integer $s= \rk( z)$ is called \emph{rank} of $z$.  For every $z \in \M$ let $\max\{z\}$ denote the largest eigenvalue of $z$, then $\max\{\cdot \}$ is a norm on $\M$ called the \emph{spectral norm}. It is worth to point out that the HSSNT $\Omega\subset\M$ associated to $\M$
is the open unit ball in $\mathcal M$ centered at the origin (with respect the spectral norm),
i.e.,
\begin{equation*}\label{Mball}
{\Omega}=\{z=\sum_{j=1}^s\lambda_j\,c_j \ |\ \max\{z\}= \max_j\{ \lmb _j\}<1\}.
\end{equation*}
The rank of $\M$ is $\rk( \M) = \max \{\rk( z)\, |\, z \in \M    \}$, moreover $\rk (\M) =\rk( \Omega) = r$.
The elements $z$ such that $\operatorname{rk}z=r$ are called \emph{regular}. If $z\in {\M}$ is regular, with spectral decomposition%
\begin{equation*}
z=\lambda_{1}\,e_{1}+\cdots+\lambda_{r}\,e_{r}\qquad(\lambda_{1}>\cdots
>\lambda_{r}>0), \label{D8}%
\end{equation*}
then $\left(  e_{1},\ldots,e_{r}\right)  $ is a \emph{(Jordan) frame of }${\M}$,
that is, a maximal sequence of pairwise orthogonal minimal tripotents.

Using the spectral decomposition, it is possible to associate to an \emph{odd} function
$f:\mathbb{R}\rightarrow \mathbb{C}$ a map  $F:\M \rightarrow \M$
as follows. Let $z\in \M$
and let
\[
z=\lambda_{1}\,c_1+\cdots+\lambda_{s}\,c_s,\quad 0<\lambda_{1}<\cdots<\lambda
_{s}
\]
be the spectral decomposition of $z$. Define the map $F$
by
\begin{equation*}\label{associatedfunction}
F\left(z\right)=f\left(\lambda_1\right)c_1+\cdots+f\left(\lambda_s\right)c_s.
\end{equation*}
If $f$ is continuous, then $F$ is continuous. If
\[
f\left(t\right)=\sum_{k=0}^{N}a_{k}\,t^{2k+1}
\]
is a polynomial, then $F$ is the map defined by
\[
F\left(z\right)=\sum_{k=0}^{N}a_{k}\,z^{\left(2k+1\right)}\qquad\left(z\in \M\right). \label{F02}
\]
If $f$ is analytic, then $F$ is real-analytic. If $f$ is given near $0$ by
\[
f\left(t\right)=\sum_{k=0}^{\infty}a_{k}\,t^{2k+1},
\]
then $F$ has the Taylor expansion near $0\in {\M}$:
\[
F\left(z\right)=\sum_{k=0}^{\infty}a_{k}\,z^{\left(2k+1\right)}. \label{F03}
\]

\begin{ex}\label{polysp}\rm
Let $P=\left(\C H^1\right)^\ell \subset \left(\C^\ell, \{,,\}\right)$ be the polydisk embedded in is its associated HPJTS $\left(\C^\ell,\{,,\}\right)$, endowed with the hyperbolic metric induced by \eqref{hypform}. Define $\tilde c_j= \left(0,\dots,\,0,\, e^{i \theta_j},\,0,\dots,\,0\right),$ $1\leq j \leq \ell$. The $\tilde{c}_j$ are mutually strongly orthogonal tripotents. Given $z=\left(\rho_1 e^{i \theta_1}, \dots,\, \rho_\ell e^{i \theta_\ell}\right) \in \C^\ell,$ $z\neq 0,$ then up to a permutation of the coordinates, we can assume $0 \leq \rho_1 \leq \rho_2 \leq \dots \leq \rho_\ell$. Let $i_1,$ $1 \leq i_1 \leq \ell,$ the first index such that $\rho_{i_1} \neq 0$ then we can write
$$
z=\rho_{i_1}\left(\tilde{c}_{i_1} + \dots + \tilde{c}_{i_2-1}\right) + \rho_{i_2}\left(\tilde{c}_{i_2} + \dots + \tilde{c}_{i_3-1}\right) + \dots + \rho_{i_s}\left(\tilde{c}_{i_s} + \dots + \tilde{c}_{i_{s+1}-1}\right)
$$
with $0 < \rho_{i_1} < \rho_{i_2} < \dots < \rho_{i_s}=\rho_\ell$ and $i_{s+1}=\ell+1$. The $c_j$'s, defined by $c_j =\tilde{c}_{i_j}+ \dots +\tilde{c}_{i_{j+1}-1}$, are still mutually strongly orthogonal tripotents and $z=\lmb_1\, c_1 + \dots + \lmb_s\, c_s$ with $\lmb_j=\rho_{i_j},$ is the spectral decomposition of $z$.

Suppose that $z=\sum_{j=1}^\ell \lmb_j\, c_j$ is a regular point, then $c_j=e^{i\theta_j}$, $1\leq j \leq l$. So the exponential map $\exp_0^P:\C^\ell\cong T_0 P\f P$ can be written as
\begin{equation}\label{expop}
\exp_0^{P} \left(z\right)=\left(\tanh \left(|z_1|\right)\frac{z_1}{|z_1|}, \dots, \tanh \left(|z_{\ell}|\right) \frac {z_\ell} {|z_\ell|}\right)= \sum^\ell_{j=1} \tanh\left(\lmb_j\right)  c_j
\end{equation}
and $\exp_0^{P}\left(0\right)=0$. Note that $\exp_0^{P} \left(z\right)= \sum^\ell_{j=1} \tanh\left(\lmb_j\right)  c_j$ is the spectral decomposition of $\exp_0^{P} \left(z\right)$. The distance from the origin of $\exp_0^{P}\left(0\right)=0$ is given by 
\begin{equation}\label{distop}
\di_\om(0,\,\exp_0^{P}(z))=\sqrt{\sum^\ell_{j=1} \lmb_j^2 }.
\end{equation}
\end{ex}

\subsection{Polar coordinates}
Let $M$ be the set ot tripotents elements of the positive Jordan triple of ${\M}$.
Then $M$ is a compact submanifold of ${\M}$ (with connected components of
different dimensions).
The \emph{height} $k$ of a tripotent element $c$ is the maximal length of a
decomposition $c=c_{1}+\cdots+c_{k}$ into a sum of pairwise orthogonal
(minimal) tripotents. Minimal tripotents have height $1$, maximal tripotents
have height $r=\operatorname{rk}{\M}$. Denote by $M_{k}$ the set of tripotents of
height $k$. \emph{If }${\M}$ \emph{is simple (that is, if }$\Omega$ \emph{is
irreducible)}, the submanifolds $M_{k}$ are the connected components of $M$.

The set $\mathcal{F}$ of frames (also called F\"{u}rstenberg-Satake boundary
of $\Omega$):
\begin{equation*}
\mathcal{F}=\left\{  \left(  c_{1},\ldots,\,c_{r}\right)  \mid c_{j}\in
M_{1},\ c_{j}\perp c_{k}\ (1\leq j<k\leq r)\right\}  , \label{P15}%
\end{equation*}
(where $c_{j}\perp c_{k}$ means orthogonality of tripotents: $T(c_{j}%
,c_{k})=0$) is a compact
submanifold of ${\M}^{r}$. The map $F: \left\{  \lambda_{1}>\cdots>\lambda_{r}>0\right\}  \times\mathcal{F}
\rightarrow {\M}_{\mathrm{reg}}\nonumber$ defined by:
\begin{equation}\label{polar coord}
\left(  \left(  \lambda_{1},\ldots,\,\lambda_{r}\right)  ,\left(  c_{1}%
,\ldots,\,c_{r}\right)  \right)  \mapsto\sum\lambda_{j}\,c_{j}%
\end{equation}
is a diffeomorphism onto the open dense set ${\M}_{\mathrm{reg}}$ of regular elements of
${\M}$, moreover its restriction%
\[
\left\{  1>\lambda_{1}>\cdots>\lambda_{r}>0\right\}  \times\mathcal{F}%
\rightarrow\Omega_{\mathrm{reg}}%
\]
is a diffeomorphism onto the set $\Omega_{\mathrm{reg}}$ of regular elements
of $\Omega$. This map plays the same role as polar coordinates in rank one.

\subsection{Geodesics and geodesic distance of HSSNT}
Let $\M$ be a HPJTS with its associated HSSNCT $\Omega\subset\M$.
For every point $z\in \M\cong T_0 \Omega$, by the polydisc Theorem (see \cite{helgason}), there exists a totally geodesic polydisc $P=\left( \C H^1\right)^r\subset \Omega$ of maximal rank $r=\rk( \Omega)$, trough the origin $0\in \M$ such that $z \in T_0 P$. Assume that $z=\sum^r_{j=1}\lmb_j\,c_j$ is a regular point then  by Proposition \ref{sub} and \eqref{expop}, the exponential map in polar coordinates is given by
\begin{equation*}%\label{expo}
\exp_0^{\Omega} \left(z\right)= \sum^r_{j=1} \tanh\left(\lmb_j\right)  c_j.
\end{equation*}
On the other hand, by \eqref{distop}, the distance from the origin of $z$ is given by
\begin{equation}\label{disto}
\di_\om(0,\,z)=\sqrt{\sum^r_{j=1} \arctanh^2\left(\lmb_j\right) }.
\end{equation}

%%%%%%%%%%%%%%%%%%%%%%%%%%%%%%%%%%%%%%%%%%%%%%%

\section{Proof of Theorem \ref{alternative}}\label{sect 3}
In order to prove Theorem \ref{alternative}  we need some properties of the symplectic duality map associated to an HSSNCT $(\Omega, g_\om)$. Recall that  the symplectic duality map  is a diffeormophism  $\Psi:\M\f\Omega\subset\M$ such that $\Psi(0)=0$ (where $0$ denotes the origin of $\M$),  $\Psi^{*}\omega_{\om}=\omega_{0}$ and hence $\Psi^{*}\frac{\omega_\om^n}{n!}=\frac{\omega_0^n}{n!}$. Therefore $\Psi$  provides  global symplectic coordinates for $(\Omega, g_{\om})$. In polar coordinates the the symplectic duality map is given by
\begin{equation*}\label{polyduality}
\Psi(t_1\,c_1+\dots+t_r\,c_r)=\left(\frac{t_1}{\sqrt{|t_1|^2+1}}\,c_1+\dots+\frac{t_r}{\sqrt{|t_r|^2+1}}\,c_r\right).
\end{equation*}
The composition of $\Psi$ with the distance \eqref{disto} is written:
\begin{equation}\label{ineqentr}
\begin{split}
\di_{\om}\left(0,\Psi(t_1v_1+\dots+t_rv_r)\right)&=\sqrt{\sum_{j=1}^r\arctanh^{2}\frac{|t_j|}{\sqrt{|t_j|^2+1}}}=\sqrt{\sum_{j=1}^r\left(\sinh^{-1}|t_j|\right)^2}.
\end{split}
\end{equation}
Let us denote $\Lambda = \{\lmb_1 >\dots>\lmb_r>0\}$. Let $F:\Lambda \times \FF \f \M _{\reg}$ the polar coordinates map given in \eqref{polar coord}. The pullback of the flat volume form $\frac{\omega_0^n}{n!}$ (see \cite[(5.1.1)]{koranyi} and \cite{unicdual}) is
\begin{equation}\label{volume form }
F^*\frac{\omega_0^n}{n!}= \Theta \wedge \prod_{j=1}^r \lmb_j^{2b+1} \prod_{1\leq j<k\leq r} (\lmb_j^2 - \lmb_k^2)^a  \ d \lmb_1 \wedge \dots \wedge d \lmb_r
\end{equation}
where $\Theta$ is a volume form on the compact manifold $\FF$. By \eqref{ineqentr} and \eqref{volume form } we get
\begin{equation}\label{the integral}
\begin{split}
&\int_{\Omega}e^{-c\,\di_{\om}(p_0,\cdot)}\frac{\omega_{\om}^n}{r!}=\int_{\M}e^{-c\,\di_\om(p_0,\Psi(x))}\,\frac{\omega_0^n(x)}{r!}\\
&\quad =K\int_{\Lambda}e^{-c\,\sqrt{ \sum_{j=1}^r\left(\sinh^{-1}\lmb_j\right)^2}} \prod_{j=1}^r \lmb_j^{2b+1} \prod_{1\leq j<k\leq r} (\lmb_j^2 - \lmb_k^2)^a  \ d \lmb_1 \wedge \dots \wedge d \lmb_r.\\
\end{split}
\end{equation}
Where $K=\int_\FF \Theta$. In order to calculate $\Entv(\Omega, g_{\hyp})$, by definition \eqref{def vol ent int}, we have to determine the minimum of $c$ such that the integral \eqref{the integral} converge. By substituting the variables $\lmb_j$ with $\sinh(t_j)$ we are reconducted to the study of the following integral:
\begin{equation*}
\begin{split}
&\int_{\Lambda}e^{-c\,\sqrt{ \sum_{j=1}^r t_j ^2}} \prod_{j=1}^r \left(\sinh t_j\right)^{2b+1} \prod_{1\leq j<k\leq r} (\sinh^2 t_j - \sinh ^2 t_k)^a\prod_{j=1}^r \left(\cosh t_j\right)\ d t_1 \wedge \dots \wedge d t_r,
\end{split}
\end{equation*}
in spherical coordinates $G:{S^{r-1}\cap\Lambda}\times \R$, $G(x,\rho)=\rho\,x$ we obtain
\begin{equation*}
\begin{split}
\int_{{S^{r-1}\cap\Lambda}}\int_{0}^\infty e^{-c\,\rho} \prod_{j=1}^r \left(\sinh \left(\rho\, x_j\right)\right)^{2\,b+1} \prod_{1\leq j<k\leq r} (\sinh^2 \left(\rho\, x_j\right) - \sinh ^2 \left(\rho\, x_k\right)^a&\cdot\\
\qquad\qquad \cdot \prod_{j=1}^r \left(\cosh \left(\rho\, x_j\right)\right)\rho^{r-1}&\ d \rho\, dv_{S^{r-1}}.\\
\end{split}
\end{equation*}
where $dv_{S^{r-1}}$ is the standard volume form on $S^{r-1}$. Since the function $e^{-\rho}$ is bounded, by dominated convergence theorem, the convergence of the previous integral is equivalent to the convergence following integral
\begin{equation*}
\begin{split}
&\int_{0}^\infty e^{-c\,\rho} \prod_{j=1}^r  e^{\rho\, (2\,b+1) \, x_j(t)}\prod_{1\leq j<k\leq r} \left(e^{2\,\rho \, x_j} - e^{2\,\rho \, x_k} \right)^a \prod_{j=1}^r e^{\rho \, x_j}\,\rho^{r-1}\  d \rho\\
&=  \int_{0}^\infty e^{-c\,\rho}   e^{\rho\, (2\,b+1) \, \sum_{j=1}^r x_j}\left(\sum_{\sigma} {\rm sgn}(\sigma)\, e^{2\,a\, \rho \, \sum_{j=1}^{r-1}\left(r-j\right)\, x_{\sigma(j)}} \right) e^{\rho \, \sum_{j=1}^r x_j}\,\rho^{r-1}\  d \rho,\\
\end{split}
\end{equation*}
where $\sigma=(\sigma_1,\dots,\sigma_r)$ is a permutation of $(1, \dots, r)$ and $x_1>\dots>x_r$. As for every $\ep >0$ and for sufficiently large $\rho$ we have $ \rho^{r-1} < e^{\ep\,\rho}$ the previous integral is convergent if and only if
\begin{equation*}
\begin{split}
c > \lim_{\ep \f 0^+}(2\,b + 2)\sum_{j=1}^r x_j + 2\,a \sum_{j=1}^{r-1}(r-j)\, x_{\sigma(j)}+\ep 
\end{split}
\end{equation*}
for every choice  of $\sigma$ and $(x_1, \dots, \,x_{r}) \in S^{r-1} \cap \Lambda$. 
Therefore the volume entropy is given by
\begin{equation*}%\label{entv hssnct}
\operatorname{Ent}\left(\Omega, g_{\om}\right)=\max_{x_1^2+\dots+x_r^2=1}\left\{2\,(b + 1)\sum_{j=1}^r  |x_j| + 2\, a \sum_{j=1}^r (r-j)\,|x_j|\right\}.
\end{equation*}
Using the method of Lagrange multipliers we see that
\begin{equation*}%\label{entv hssnct1}
\max_{x_1^2+\dots+x_r^2=1}\left\{2\,(b + 1)\sum_{j=1}^r |x_j| + 2\, a \sum_{j=1}^r (r-j)\,|x_j|\right\}= { 2\sqrt{\sum _{j=1}^{r}
 \left( b+1+a \left( r-j \right)  \right) ^{2}}}
 \end{equation*}
attained for 
\begin{equation*}%\label{entv hssnct2}
x_k= \frac {b+1+a \left( r-k \right) }{ \sqrt{\sum _{j=1}^{r}
 \left( b+1+a \left( r-j \right)  \right) ^{2}}}, \quad k=1,\dots,r.
\end{equation*}
and this concludes the proof of Theorem \ref{alternative}.

\section{proof of Theorem \ref{thm prod}}\label{sec 4}
Let now $\Omega=\prod_{k=1}^\ell \Omega_k$ be a product of irreducible bounded symmetric domains and let $(a_k,b_k,r_k)$ be the invariants associated to $\Omega_k$. As the symplectic duality map of $\Omega$ is the product of the symplectic duality map of each $\Omega_k$, arguing as in the proof of Theorem \ref{alternative}, we are reconducted to study of the converging of the following integral 
\begin{equation*}
\begin{split}
&\int_{0}^\infty e^{-c\,\rho}   e^{\rho\, \sum_{k=1}^\ell \left( (2\,b_k+2) \, \sum_{j=1}^r x_{k,j}\right)}\prod_{k=1}^\ell \left(\sum_{\sigma_k}  {\rm sgn}(\sigma_k)\,  e^{2\,a_k\, \rho \, \sum_{j=1}^{r_k-1}\left(r_k-j\right)\, x_{k,\sigma_k(j)}} \right)\rho^{r-1}\  d \rho,
\end{split}
\end{equation*}
where $\sigma_k=(\sigma_k(1),\dots,\sigma_k(r_k))$ is a permutation of $(1, \dots, r_k)$ and $x_{k,1}>\dots>x_{{k,r_k}}$.
Therefore the volume entropy is given by
\begin{equation*}%\label{entv hssnct}
\operatorname{Ent}\left(\Omega, g_{\om}\right)=\max_{\sum_{k,j} x_{k,j}^2=1}\left\{2\,\sum_{k=1}^\ell \left( (\,b_k+1) \, \sum_{j=1}^{r_k} |x_{k,j}| + \, a_k \sum_{j=1}^{r_k} (r_k-j)\,|x_{k,j}| \right)\right\}.
\end{equation*}
Using the method of Lagrange multipliers we see that
\begin{equation*}%\label{entv hssnct1}
\max_{\sum_{k,j} x_{k,j}^2=1}\left\{2\,\sum_{k=1}^\ell \left( (\,b_k+1) \, \sum_{j=1}^{r_k} |x_{k,j}| + \, a_k \sum_{j=1}^{r_k} (r_k-j)\,|x_{k,j}| \right)\right\}
\end{equation*}
\begin{equation*}%\label{entv hssnct1}
=2{ \sqrt{\sum_{k=1}^\ell \sum _{j=1}^{r_k}
 \left( b_k+1+a_k \left( r_k-j \right)  \right) ^{2}}}
\end{equation*}
\begin{equation*}%\label{entv hssnct1}
=\sqrt{\sum_{k=1}^\ell\Entv^2(\Omega_k,g_k)}
\end{equation*}

attained for 
\begin{equation*}%\label{entv hssnct2}
x_{k,j}= \frac {b_k+1+a_k \left( r_k-j \right) }{ \sqrt{\sum_{k=1}^\ell \sum _{j=1}^{r_k}
 \left( b_k+1+a_k \left( r_k-j \right)  \right) ^{2}}}, \quad k=1,\dots,\ell,\ j=1,\dots,r_k.
\end{equation*}
and this concludes the proof of Theorem \ref{thm prod}.

\begin{rmk}\rm
The volume entropy $\Entv\left( \Omega, g_B \right)$ of a irreducible HSSNCT with respect to the Bergman metric $g_B=\gamma\, g_\om$ is given by
$$\Entv\left( \Omega, g_B \right)=\frac{\Entv\left( \Omega, g_\om \right)}{\sqrt\gamma}$$
This follows by the definition of the volume entropy \eqref{def vol ent int} and by the fact that the distance with respect the Bergman metric satisfies $\di_B\left(p,q\right)={\sqrt\gamma}\,\di_\om\left(p,q\right)$.
\end{rmk}

\end{document}